\documentclass[10pt,a4paper]{article}
\usepackage[T1]{fontenc}
\usepackage{graphicx}
\usepackage{amssymb}
\usepackage{amsthm}
\usepackage{amsmath}
\usepackage{hyperref}
\usepackage{array}
\usepackage{booktabs}
\usepackage{float}
\usepackage{algorithm}
\usepackage[noend]{algpseudocode}
\PassOptionsToPackage{numbers}{natbib}
\usepackage[preprint]{neurips_2024}
\newtheorem{theorem}{Theorem}[section]

\author{%
	Xu Duan$^1$, Dongmei Chen$^{1}$\\
	$^1$The University of Texas at Austin\\
	\texttt{\{xuduan, dmchen\}@utexas.edu}}
\title{Trajectory-Optimized Density Control with Flow Matching}
\begin{document}
	\maketitle

	\begin{abstract}
		Optimal transport (OT) and Schrödinger bridge (SB) problems have emerged as powerful frameworks for transferring probability distributions with minimal cost. However, existing approaches typically focus on endpoint matching while neglecting critical path-dependent properties—particularly collision avoidance in multi-agent systems—which limits their practical applicability in robotics, economics, and other domains where inter-agent interactions are essential. Moreover, traditional density control methods often rely on independence assumptions that fail to capture swarm dynamics. We propose a novel framework that addresses these limitations by employing flow matching as the core modeling tool, where the flow model co-evolves with the control policy. Unlike prior methods that treat transport trajectories as mere interpolations between source and target distributions, our approach explicitly optimizes over the entire transport path, enabling the incorporation of trajectory-dependent costs and collision avoidance constraints. Building on recent advances in neural PDE solvers and conditional generative modeling, we develop a simulation-free method that learns transport policies without requiring explicit simulation of agent interactions. Our framework bridges optimal transport theory with mean-field control, providing a principled approach to multi-agent coordination problems where both endpoint alignment and path properties are critical. Experimental results demonstrate that our method successfully generates collision-free transport plans while maintaining computational efficiency comparable to standard flow matching approaches.
	\end{abstract}
	
	\section{Introduction}
	Optimal transport (OT) is a mathematical framework for transferring one probability distribution to another with minimal cost. Classical OT primarily focuses on the source and target distributions, with the cost function determined by their alignment. However, in many applications, the trajectory of the transport itself significantly contributes to the overall loss. This extended perspective is particularly important in domains such as economics \citep{Achdou2014}, opinion dynamics \citep{Schweighofer2020}, and robotics \citep{Liu2018}, where the path of transport directly influences outcomes.
	
	This class of problems is closely related to Schrödinger Bridge (SB) theory, an entropy-regularized variant of OT \citep{Chen2020,Lonard2013}. SB formulates the problem as a minimum-effort control task between two prescribed distributions over a finite time horizon. When dynamics are deterministic, OT provides an appropriate framework, whereas SB is more suitable in stochastic settings \citep{Chen2016}. A key limitation of existing density control methods, however, lies in their reliance on independence assumptions across individuals, particularly in swarm-control scenarios. As a result, critical swarm properties—most notably collision avoidance—are often left unmodeled \citep{Krishnan2018,Inoue2021,Elamvazhuthi2019}, which undermines their applicability in practice.
	
	A related line of work arises from mean-field games (MFG) \citep{Bensoussan2013}, which explicitly model interactions in large populations of agents. Recent developments have examined these problems through the lens of conditioning or inverse problems in diffusion or flow-based models \citep{Daras2024}. For example, \citet{Kerrigan2024} propose a simulation-free, flow-based method for conditional generative modeling using conditional optimal transport (COT). Their approach couples an arbitrary source distribution to a specified target via a triangular COT plan, approximating the geodesic path of measures to construct a conditional generative model. Similarly, \citet{Whang2021} design an approximate inference scheme that expresses the conditional target as a composition of two flow models, thereby enabling stable variational inference training without adversarial techniques. Yet, as with SB and swarm-control approaches, these methods do not explicitly encode collision avoidance, limiting their utility in multi-agent settings where inter-agent interactions are critical.
	
	The rapid advancement of machine learning has also spurred progress in numerical methods for these problems. Seminal works such as \citet{Ruthotto2020} approximate solutions with deep neural networks (DNNs) by directly penalizing PDE violations. \citet{Mei2025a} employ a time-reversal method to formulate stochastic optimal control in a forward–backward SDE (FBSDE) framework. A deterministic variant is considered by \citet{Zhou2025b}, who leverage entropy differentials—namely the score function—to derive deterministic forward–backward characteristics for mean-field control systems.
	
	Our work is strongly motivated by \citet{Liu2022} but departs in two key aspects:
	(1) we employ flow matching as the core modeling tool, and
	(2) our flow matching model co-evolves with the policy.
	The ultimate goal of our framework is to identify transport paths that minimize cost not only at the endpoints, but also along the trajectory itself.

	\section{Background}
	\subsection{Preliminary on optimal transport}
	Given two nonnegative measures $\mu, \nu$ on $\mathbb{R}^d$ having equal total mass (often assumed to be probability distributions), the Monge's formulation of optimal transport seeks a transport map
	\begin{equation}
		T : \mathbb{R}^d \to \mathbb{R}^d : x \mapsto T(x)
	\end{equation}
	
	from $\mu$ to $\nu$ in the sense $T_\# \mu = \nu$, that incurs minimum cost of transportation
	\begin{equation}\label{monge}
	    \int c(x, T(x)) \mu(dx).
	\end{equation}
	
	Here, $c(x, y)$ stands for the transportation cost per unit mass from point $x$ to $y$. The dependence of the total transportation cost on $T$ is highly nonlinear, complicating early analyses to the problem [10]. This problem was later relaxed by Kantorovich, where, instead of a transport map, a joint distribution $\pi$ on the product space $\mathbb{R}^d \times \mathbb{R}^d$ is sought. Let $\Pi(\mu, \nu)$ be the set of joint distributions of $\mu$ and $\nu$, then the Kantorovich formulation of OT reads
	\begin{equation}\label{kant}
		\inf_{\pi \in \Pi(\mu, \nu)} \int_{\mathbb{R}^d \times \mathbb{R}^d} c(x, y) \pi(dx dy).
	\end{equation}
	
	Both the Monge's and the Kantorovich's formulations are "static" focusing on "what goes where". It turns out that the OT problem can also be cast as a dynamical, temporal dimension. In particular, when $c(x, y) = \frac{1}{2} \|x - y\|^2$, OT can be formulated as a stochastic control problem
	\begin{equation}\label{dyn}
		\inf_{v \in \mathcal{V}} \mathbb{E} \left\{ \int_0^1 \frac{1}{2} \|v(t, x^v(t))\|^2 dt \right\},
	\end{equation}
	where
	\begin{equation}
		x^v(t) = x(0) + \int_0^t v(s, x^v(s)) ds,
	\end{equation}
	\begin{equation}
		x^v(0) \sim \mu, \quad x^v(1) \sim \nu.
	\end{equation}
	
	Here $\mathcal{V}$ represents the family of admissible state feedback control laws. Note that this control problem differs from standard ones \ref{monge} in that the terminal constraint $x^v(1) \sim \nu$, meaning $x^v(1)$ follows distribution $\nu$, is unconventional. In \ref{monge}, the goal is to find an optimal control policy to drive the system from an uncertain initial state $x^v(0) \sim \mu$ to an uncertain target state $x^v(1) \sim \nu$. The solution to \ref{dyn} specifies how to move mass over time from configuration $\mu$ to $\nu$, providing more resolution to the optimal transport plan.
	\subsection{Preliminary on schrödinger bridge}
	In 1931/32, Schrödinger posed the following problem \cite{Schrdinger1931, Schrdinger1932}: A large number $N$ of independent Brownian particles in $\mathbb{R}^d$ is observed to have an empirical distribution approximately equal to $\mu$ at time $t = 0$, and at some later time $t = 1$ an empirical distribution approximately equal to $\nu$. Suppose that $\nu$ differs from what it should be according to the law of large numbers, namely
	\begin{equation}
	\int q_e(0, x; 1, y) \mu(dx),
	\end{equation}
	where
	\begin{equation}
	q_e(s, x; t, y) = (2\pi)^{-d/2} |t - s|^{-d/2} \exp\left(-\frac{\|x - y\|^2}{2(t - s)}\right)
	\end{equation}
	denotes the scaled Brownian transition probability density. It is apparent that the particles have been transported in an unlikely way. But of the many unlikely ways in which this could have happened, which one is the most likely?
	
	This problem can be understood in the modern language of large deviation theory as a problem [34] of determining a probability law $P$ on the path space $\Omega = C([0, 1], \mathbb{R}^d)$ that minimizes the relative entropy (a.k.a., Kullback-Leibler divergence)
	\begin{equation}
	KL(P \| Q) := \int_\Omega \log\left(\frac{dP}{dQ}\right) dP.
	\end{equation}
	Here $Q$ is the probability law induced by the Brownian motion and $P$ is chosen among probability laws that are absolutely continuous with respect to $Q$ and have the prescribed marginals. The solution to this optimization problem is referred to as the Schrödinger bridge. Existence and uniqueness of the minimizer have been proven in various degrees of generality by Fortet \cite{Fortet1940}, Beurling \cite{Beurling1960}, Föllmer \cite{Fllmer1988}.
	\section{Theory}
	\subsection{Generalized SB-FBSDEs}
	
	According to \cite{Liu2022}, the Generalized FBSDEs is as follows  
	
	\begin{theorem}[Generalized SB-FBSDEs]
		Suppose $\Psi, \hat{\Psi} \in C^{2,1}$ and let $f, F$ satisfy usual growth and Lipschitz conditions. Consider the following nonlinear FK transformations applied to (9):
		\begin{align}
			Y_t &\equiv Y(X_t, t) = \log \Psi(X_t, t), \quad Z_t \equiv Z(X_t, t) = \nabla \log \Psi(X_t, t), \\
			\hat{Y}_t &\equiv \hat{Y}(X_t, t) = \log \hat{\Psi}(X_t, t), \quad \hat{Z}_t \equiv \hat{Z}(X_t, t) = \sigma \nabla \log \hat{\Psi}(X_t, t).
		\end{align}
		where $X_t$ follows (3a) with $X_0 \sim \rho_0$. Then, the resulting FBSDEs system takes the form:
		\begin{equation}\label{FBSDE}
			\text{FBSDEs w.r.t. (3a)} \quad 
			\begin{cases}
				dX_t = (f_t + \sigma Z_t)dt + \sigma dW_t \\
				dY_t = \Big(\tfrac{1}{2}\|Z_t\|^2 + F_t \Big)dt + Z_t^T dW_t \\
				d\hat{Y}_t = \Big(\tfrac{1}{2}\|\hat{Z}_t\|^2 + \nabla \cdot (\sigma \hat{Z}_t - f_t) + 2 \hat{Z}_t^T Z_t - F_t \Big)dt + \hat{Z}_t^T dW_t
			\end{cases}
		\end{equation}
		
		Now, consider a similar transformation in (9) but instead w.r.t. the ``reversed'' SDE $X_s \sim (3b)$ and $X_0 \sim \rho_{\text{target}}$, i.e. $Y_s \equiv Y(X_s, s) = \log \Psi(X_s, s)$, and etc. The resulting FBSDEs system reads:
		\begin{equation}
			\text{FBSDEs w.r.t. (3b)} \quad 
			\begin{cases}
				dX_s = (-f_s + \sigma \hat{Z}_s) ds + \sigma dW_s \\
				dY_s = \Big(\tfrac{1}{2}\|Z_s\|^2 + \nabla \cdot (\sigma Z_s + f_s) + 2 Z_s^T \hat{Z}_s - F_s \Big) ds + Z_s^T dW_s \\
				d\hat{Y}_s = \Big(\tfrac{1}{2}\|\hat{Z}_s\|^2 + F_s \Big) ds + \hat{Z}_s^T dW_s
			\end{cases}
		\end{equation}
		
		Since $Y_t + \hat{Y}_t = \log \rho(X, t)$ by construction, the functions $f_t$ and $F_t$ in (11) take the arguments
		\[
		f_t := f_t(X_t, \exp(Y_t + \hat{Y}_t)), \qquad F_t := F_t(X_t, \exp(Y_t + \hat{Y}_t)).
		\]
		Similarly, we have $f_s := f_s(X_s, \exp(Y_s + \hat{Y}_s))$ and $F_s := F_s(X_s, \exp(Y_s + \hat{Y}_s))$ in (12).
	\end{theorem}
	
	%\begin{proof}
	%	See Appendix A.3.3.
	%\end{proof}
	which is a generalization of \cite{Chen2023} by introducing nontrivial MF interaction $F$.  The corresponding Iterative Proportional Fitting (IPF) loss function is defined as
	\begin{align}
		\mathcal{L}_{\text{IPF}}(\theta) &= \int_0^T \mathbb{E} \left[ \frac{1}{2} \|\mathcal{Z}_\theta(\bar{X}_{s,s})\|^2 + \mathcal{Z}_\theta(\bar{X}_{s,s})^\top \mathcal{Z}_\phi(\bar{X}_{s,s}) + \nabla \cdot (\sigma \mathcal{Z}_\theta(\bar{X}_{s,s}) + f) \right] ds, \label{eq:6a} \\
		\mathcal{L}_{\text{IPF}}(\phi) &= \int_0^T \mathbb{E} \left[ \frac{1}{2} \|\mathcal{Z}_\phi(\bar{X}_{t,t})\|^2 + \mathcal{Z}_\phi(\bar{X}_{t,t})^\top \mathcal{Z}_\theta(\bar{X}_{t,t}) + \nabla \cdot (\sigma \mathcal{Z}_\phi(\bar{X}_{t,t}) - f) \right] dt. \label{eq:6b}
	\end{align}
	
	\subsection{Temporal Difference objective $\mathcal{L}_{TD}$}

	Discretizing \ref{FBSDE}  with some fixed step size $\delta t$ yields
	\begin{equation}
		Y_{t+\delta t}^\theta = Y_t^\theta + \left( \frac{1}{2} \| Z_t^\theta \|^2 + F_t \right) \delta t + Z_t^\theta S W_t, \quad \delta W_t \sim \mathcal{N}(0, \delta t I)
	\end{equation}
	which resembles a (non-discounted) Temporal Difference (TD) \cite{Todorov2009, Lutter2021} except that, in addition to standard ''reward'' (in terms of control and state costs), we also have a stochastic term. This stochastic term, which vanishes in the vanilla Bellman equation upon taking expectations, plays a crucial role in characterizing the inherited stochasticity of the value function $Y_t$. With this interpretation in mind, we can construct suitable TD targets for our FBSDEs systems as shown below.
	
	\subsection*{Proposition 3 (TD objectives). The single-step TD targets take the form:}
	\begin{align}
		\widehat{\text{TD}}_{t+\delta t}^{\text{single}} &= \widehat{Y}_t^\phi + \left( \frac{1}{2} \| \widehat{Z}_t^\phi \|^2 + \nabla \cdot (\sigma \widehat{Z}_t^\phi - f_t) + \widehat{Z}_t^\phi \text{T} \widehat{Z}_t^\phi - F_t \right) \delta t + \widehat{Z}_t^\phi \text{T} \delta W_t, \\
		\text{TD}_{s+\delta s}^{\text{single}} &= Y_s^0 + \left( \frac{1}{2} \| Z_s^0 \|^2 + \nabla \cdot (\sigma Z_s^0 + f_s) + Z_s^0 \text{T} \widehat{Z}_s^\phi - F_s \right) \delta s + Z_s^0 \text{T} \delta W_s,
	\end{align}
	with $\widehat{\text{TD}}_0 := \log \rho_0 - Y_0^\theta$ and $\text{TD}_0 := \log \rho_{\text{target}} - \widehat{Y}_0^\phi$, and the multi-step TD targets take the forms:
	\begin{align}
		\widehat{\text{TD}}_{t+\delta t}^{\text{multi}} &= \widehat{\text{TD}}_0 + \sum_{\tau=0}^{t} \delta \widehat{Y}_\tau, \\
		\text{TD}_{s+\delta s}^{\text{multi}} &= \text{TD}_0 + \sum_{\tau=0}^{s} \delta Y_\tau,
	\end{align}
	where $\delta \widehat{Y}_\tau := \widehat{\text{TD}}_{t+\delta t}^{\text{single}} - \widehat{Y}_t$ and $\delta Y_s := \text{TD}_{s+\delta s}^{\text{single}} - Y_s$. Given these TD targets, we can construct
	\begin{align}
		\mathcal{L}_{TD}(\theta) &= \sum_{s=0}^{T} \mathbb{E} \left[ \left| Y_\theta(\widehat{X}_s, s) - \text{TD}_s \right| \right] \delta s, \\
		\mathcal{L}_{TD}(\phi) &= \sum_{t=0}^{T} \mathbb{E} \left[ \left| \widehat{Y}_\phi(X_t, t) - \widehat{\text{TD}}_t \right| \right] \delta t.
	\end{align}
	
	\subsection{Flow matching objective $\mathcal{L}_{FM}$}
	The flow matching objective aims to train a neural velocity field $u_{\theta}(t, x)$ to approximate the true velocity field $u(t, x)$ that governs the transformation between two probability distributions $p_0$ and $p_1$. Formally, the objective is defined as
	\begin{equation}
		\mathcal{L}_{\mathrm{FM}}(\theta)
		\;=\;
		\mathbb{E}_{t \sim \mathcal{U}[0,1]}
		\;\mathbb{E}_{x \sim p_t}\big[\,\|u_{\theta}(t, x) - u(t, x)\|_2^2\,\big],
	\end{equation}
	where $p_t$ denotes the intermediate distribution at time $t \in [0, 1]$, and $\mathcal{U}[0,1]$ is the uniform distribution over the time interval.
	
	In the simplest linear interpolation setting, the samples evolve along straight-line trajectories connecting the source and target distributions:
	\begin{equation}
		x = (1 - t)X_0 + tX_1,
	\end{equation}
	where $X_0 \sim p_0$ and $X_1 \sim p_1$. The corresponding ground-truth velocity field is given by
	\begin{equation}
		u(t, x)
		\;=\;
		\mathbb{E}\big[X_1 - X_0 \,\big|\, X_t = x\big],
	\end{equation}
	which represents the expected displacement of particles conditioned on their location $x$ at time $t$. Intuitively, this formulation encourages the learned flow $u_{\theta}$ to match the true transport dynamics along the interpolation path between $p_0$ and $p_1$.

	\subsection{Algorithm}
	The algorithm performs gradient descent on forward and backward trajectories alternately, as detailed in Algorithm \ref{algo}.
	\begin{algorithm}[H]\label{algo}
		\caption{Learning Algorithm}
		\begin{algorithmic}[1]
			\Require $(\hat{Y}_0, \hat{Y}_\phi, \sigma \nabla \hat{Y}_\theta, \sigma \nabla \hat{Y}_\phi)$ for critic or $(Y_\theta, Y_\phi, Z_\theta, Z_\phi)$ for actor-critic parametrization.
			\Repeat
			\State Sample $\mathbf{X}^\theta = \{ \mathbf{X}_t^\theta, \mathbf{Z}_t^\theta, \delta \mathbf{W}_t \}_{t \in [0, T]}$ from the forward SDE (11a); add $\mathbf{X}^\theta$ to replay buffer $\mathcal{B}$.
			\For{$k = 1$ to $K$ do}
			\State Sample on-policy $\mathbf{X}^\theta_{\text{on}}$ and off-policy $\mathbf{X}^\theta_{\text{off}}$ samples respectively from $\mathbf{X}^\theta$ and $\mathcal{B}$.
			\State Compute $\mathcal{L}(\phi) = \mathcal{L}_{\text{IPF}}(\phi; \mathbf{X}^\theta_{\text{on}}) + \mathcal{L}_{\text{TD}}(\phi; \mathbf{X}^\theta_{\text{off}}) + \mathcal{L}_{\text{FM}}(\phi; \mathbf{X}^\theta_{\text{on}})$.
			\State Update $\phi$ with the gradient $\nabla_\phi \mathcal{L}(\phi)$.
			\EndFor
			\State Sample $\mathbf{\hat{X}}^\phi = \{ \mathbf{\hat{X}}_s^\phi, \mathbf{\hat{Z}}_s^\phi, \delta \mathbf{W}_s \}_{s \in [0, T]}$ from the backward SDE (12a); add $\mathbf{\hat{X}}^\phi$ to replay buffer $\hat{\mathcal{B}}$.
			\For{$k = 1$ to $K$ do}
			\State Sample on-policy $\mathbf{\hat{X}}^\phi_{\text{on}}$ and off-policy $\mathbf{\hat{X}}^\phi_{\text{off}}$ samples respectively from $\mathbf{\hat{X}}^\phi$ and $\hat{\mathcal{B}}$.
			\State Compute $\theta(\theta) = \mathcal{L}_{\text{IPF}}(\theta; \mathbf{\hat{X}}^\phi_{\text{on}}) + \mathcal{L}_{\text{TD}}(\theta; \mathbf{\hat{X}}^\phi_{\text{off}}) + \mathcal{L}_{\text{FM}}(\theta; \mathbf{\hat{X}}^\phi_{\text{on}})$.
			\State Update $\theta$ with the gradient $\nabla_\theta \mathcal{L}(\theta)$.
			\EndFor
			\Until{converges}
		\end{algorithmic}
	\end{algorithm}
	\section{Numerical results}
	We validate our method on crowd navigation problems. Specifically, we considers problems from \cite{Ruthotto2020} and \cite{Lin2021}. 
	\subsection{Generalized Mixture of Gaussian}
	We consider a transport problem where agents must navigate from a centralized initial distribution to a target distribution arranged in a circular pattern while avoiding obstacles. The source distribution $p_0$ is specified as a standard bivariate Gaussian centered at the origin:
	\begin{equation}
		p_0 = \mathcal{N}(\mathbf{0}, \mathbf{I}_2),
	\end{equation}
	where $\mathbf{I}_2$ denotes the $2 \times 2$ identity matrix.
	
	The target distribution $p_T$ is a Gaussian mixture model (GMM) with $N=8$ equally-weighted components arranged uniformly on a circle of radius $r=16$. Each component $i \in \{1,\ldots,8\}$ is a bivariate Gaussian with unit variance positioned at angle $\theta_i = 2\pi i/8$:
	\begin{equation}
		p_T = \frac{1}{8}\sum_{i=1}^{8} \mathcal{N}\left(\boldsymbol{\mu}_i, \mathbf{I}_2\right),
	\end{equation}
	where $\boldsymbol{\mu}_i = r[\cos(\theta_i), \sin(\theta_i)]^\top$.
	
\begin{figure}[H]
		\centering
		\includegraphics[scale=0.4]{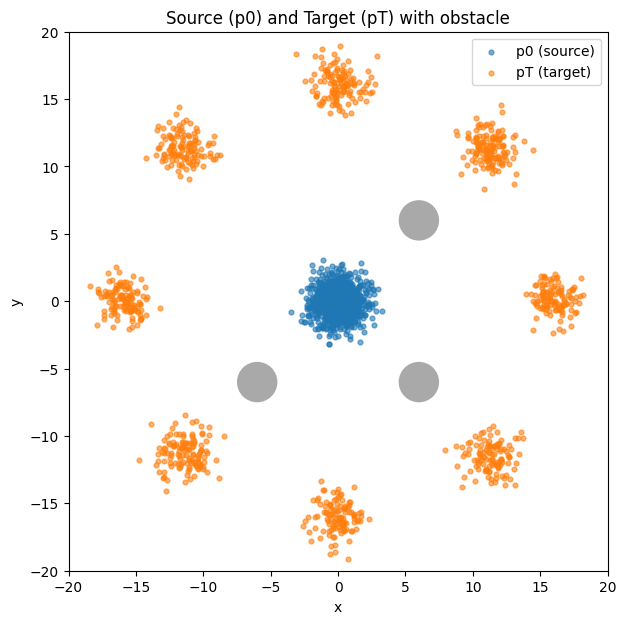}
		\caption{Problem setup for the Gaussian mixture model transport task. The source distribution (center) is a standard Gaussian centered at the origin, while the target distribution consists of eight equally-weighted Gaussian components arranged uniformly on a circle of radius $16$. Three circular obstacles of radius $1.5$ are strategically positioned to obstruct direct transport paths.}
		\label{fig:gmm_setup}
	\end{figure}

	To test collision avoidance capabilities, we introduce three circular obstacles of radius $1.5$ centered at positions $(6,6)$, $(6,-6)$, and $(-6,-6)$. These obstacles are strategically positioned along natural transport paths between the source and target distributions, requiring agents to coordinate their movements and navigate around constraints while achieving the desired distributional transfer. This configuration presents a challenging multi-agent coordination problem in which naive transport plans would result in collisions. Figure~\ref{fig:gmm_setup} illustrates the complete problem setup, showing the source distribution, the eight target GMM components, and the obstacle configuration.
	
\begin{figure}[H]
		\centering
		\includegraphics[scale=0.25]{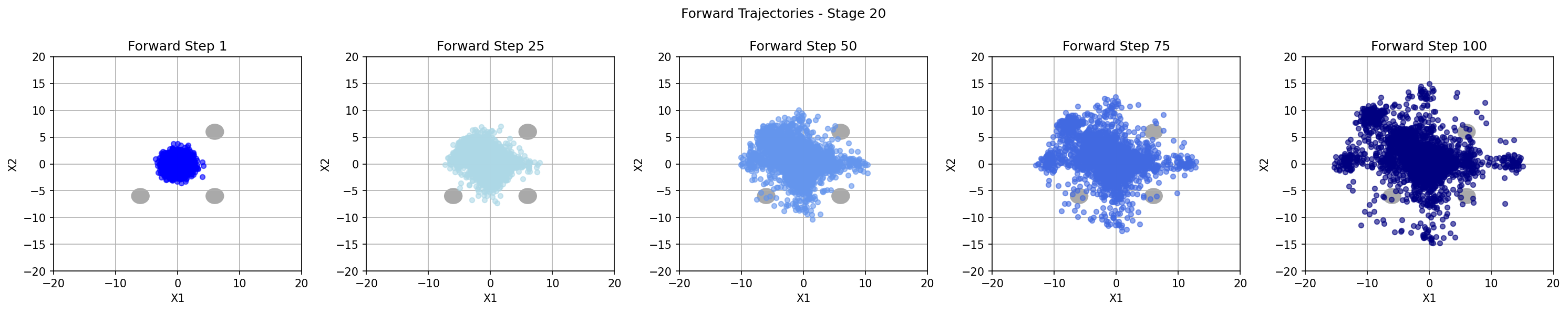}
		\includegraphics[scale=0.25]{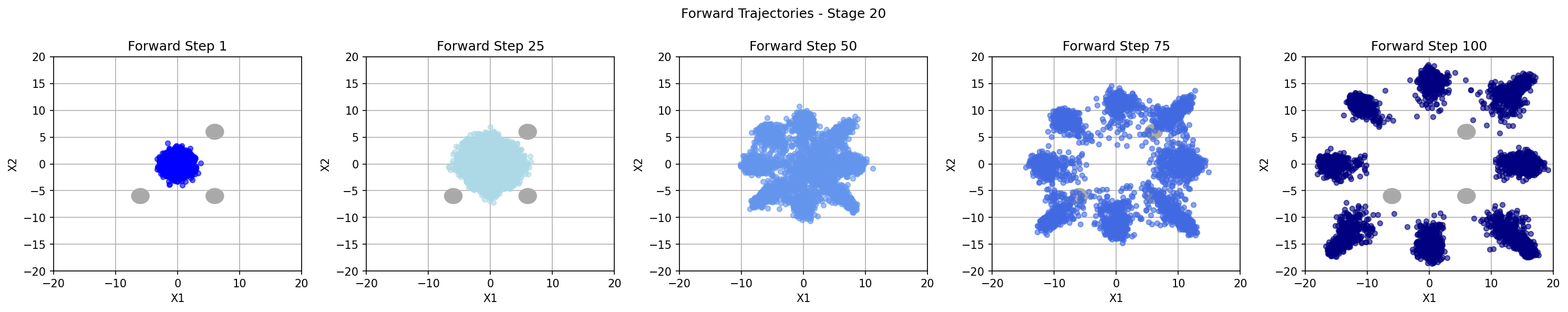}
		\caption{Forward transport for the GMM problem at stage 20. Top row: Results from the method in \cite{Liu2022}. Bottom row: Results from the proposed method. Left column: Agent trajectories navigating from the source to target while avoiding obstacles. Right column: Final distribution of agents, demonstrating successful coverage of all eight target components.}
		\label{fig:gmm_forward}
	\end{figure}

\begin{figure}[H]
		\centering
		\includegraphics[scale=0.25]{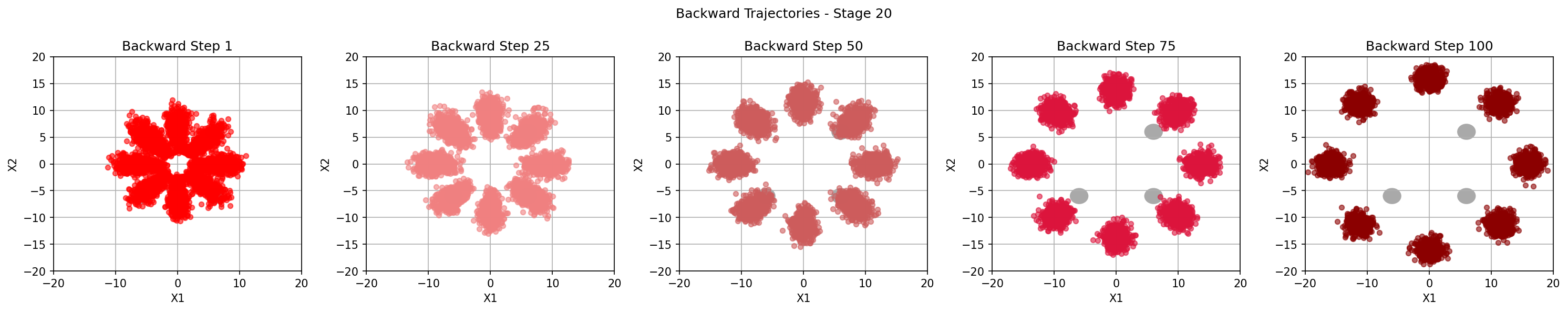}
		\includegraphics[scale=0.25]{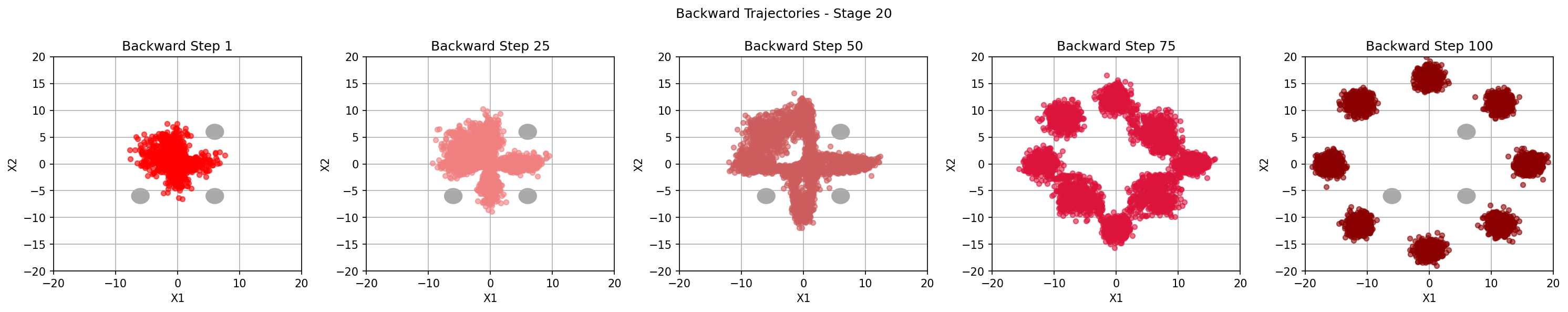}
		\caption{Backward transport for the GMM problem at stage 20. Top row: Results from the method in \cite{Liu2022}. Bottom row: Results from the proposed method. Left column: Agent trajectories from the distributed target back to the source. Right column: Final distribution showing successful reconcentration at the origin.}
		\label{fig:gmm_backward}
	\end{figure}

	The forward transport results are presented in Figure~\ref{fig:gmm_forward}, which compares agent trajectories and final distributions between the method from \cite{Liu2022} (top row) and the proposed approach (bottom row). Our method demonstrates successful navigation around obstacles while achieving uniform coverage of all eight target components. The corresponding backward transport results are shown in Figure~\ref{fig:gmm_backward}, illustrating the reverse transport from the distributed target back to the concentrated source distribution.
	
	\subsection{V-Neck Corridor}
	We consider a narrow corridor navigation problem where agents must traverse from one side to the other through a constrained passage. The source distribution $p_0$ is a concentrated Gaussian located at $(-7, 0)$:
	\begin{equation}
		p_0 = \mathcal{N}\left([-7, 0]^\top, 0.2\mathbf{I}_2\right),
	\end{equation}
	and the target distribution $p_T$ is symmetrically positioned at $(7, 0)$:
	\begin{equation}
		p_T = \mathcal{N}\left([7, 0]^\top, 0.2\mathbf{I}_2\right).
	\end{equation}
	
\begin{figure}[H]
		\centering
		\includegraphics[scale=0.4]{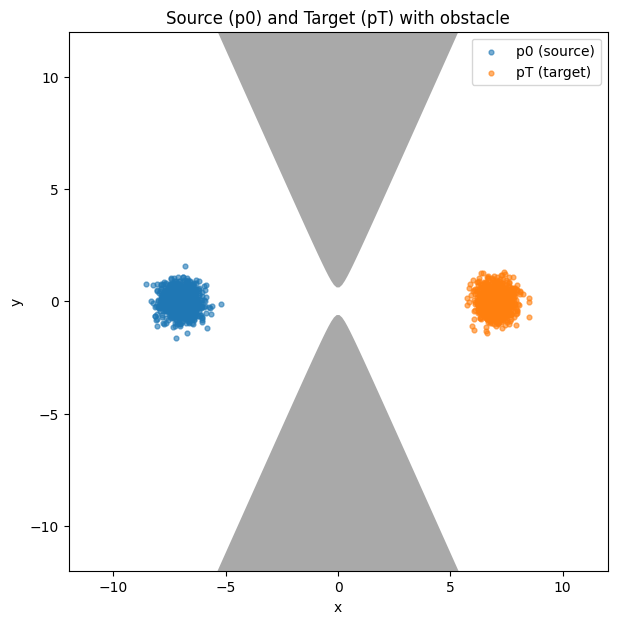}
		\caption{Problem setup for the V-neck corridor navigation task. The source distribution (left) and target distribution (right) are connected by a V-shaped constrained passage, creating a bottleneck through which all agents must coordinate their movement.}
		\label{fig:vneck_setup}
	\end{figure}

	The obstacle configuration forms a V-shaped corridor or bottleneck that agents must navigate through. The constraint is parameterized by $c^2 = 0.36$ and a coefficient $\alpha = 5$, which define the geometry of the admissible region. This setup creates a challenging scenario where agents starting from a concentrated distribution must coordinate to pass through a narrow passage without collision, then re-concentrate at the target location. Unlike the GMM scenario where agents can disperse to multiple targets, this problem requires all agents to funnel through the same constrained region, testing the method's ability to handle high-density interactions and maintain collision avoidance under severe spatial constraints. The V-neck corridor geometry is illustrated in Figure~\ref{fig:vneck_setup}.
	
\begin{figure}[H]
		\centering
		\includegraphics[scale=0.25]{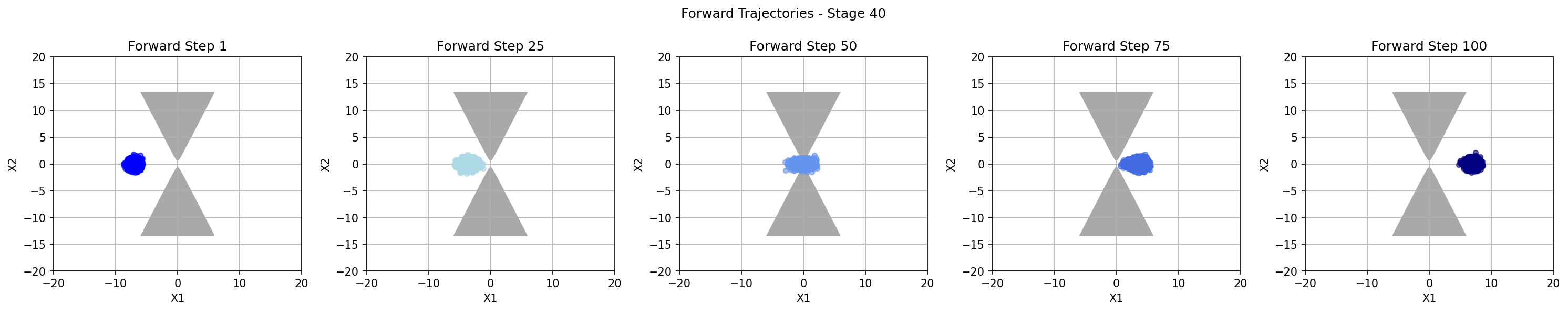}
		\includegraphics[scale=0.25]{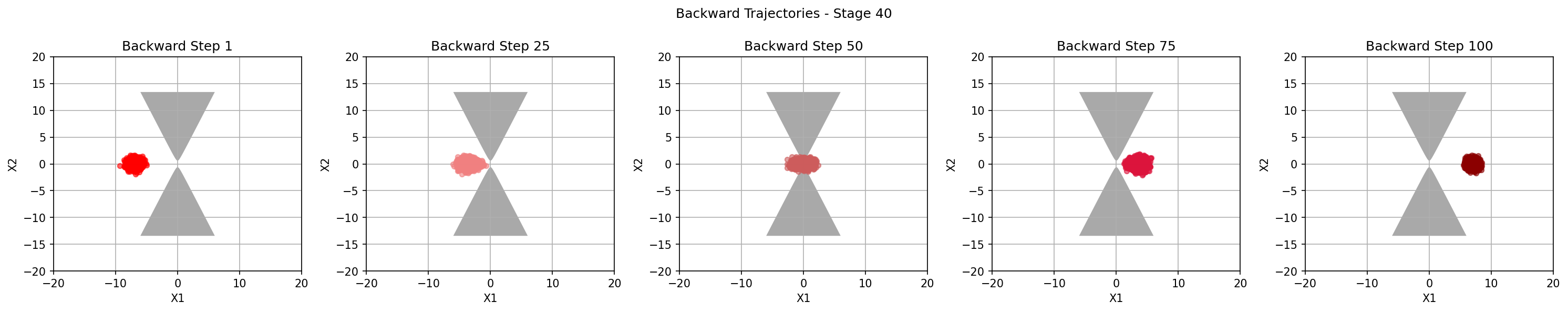}
		\caption{Transport results for the V-neck corridor problem at stage 20. Left: Forward transport showing agents successfully navigating through the bottleneck from left to right. Right: Backward transport demonstrating successful passage in the reverse direction.}
		\label{fig:vneck_transport}
	\end{figure}

	Figure~\ref{fig:vneck_transport} shows the forward and backward transport results, demonstrating successful navigation through the narrow corridor in both directions.
	
	\bibliography{export}
	\bibliographystyle{plainnat}
	\section{Appendix: Experimental Details}
	Both the flow matching network and policy network use a two-layer neural network with SiLu activation and a hidden size of 64. The hyperparameters are consistent across both experimental scenarios: the state space dimension is $d=2$, the noise level is $\sigma=1$, and the time horizon is $T=1$. We employ a discrete time step of $\delta t = 0.01$, resulting in 100 diffusion steps. For each alternating optimization stage, we use $K=250$ sample points. The complete training procedure runs for 20 alternating stages, with each stage involving 1,000 training steps, yielding a total of 20,000 training steps for both the GMM and V-neck corridor experiments.
	
\end{document}